\documentclass[12pt]{amsart}

\usepackage{amssymb}
\usepackage{amsmath}
\usepackage{amsthm}
\usepackage{amsfonts}

\usepackage{a4wide}
\usepackage{longtable}
\usepackage{multirow}
\usepackage{xcolor}

\newtheorem{theorem}{Theorem}
\newtheorem{corollary}{Corollary}

\theoremstyle{definition}
\newtheorem{definition}{Definition}
\newtheorem{example}{Example}

\newtheorem{hypothesis}{Hypothesis}

\newtheorem{warning}{Warning}

\numberwithin{equation}{section}

\newcommand{\nocontentsline}[3]{}
\newcommand{\tocless}[2]{\bgroup\let\addcontentsline=\nocontentsline#1{#2}\egroup}

\usepackage{hyperref}
\hypersetup{
    colorlinks=true,
     citecolor=blue,      
     urlcolor=blue,          
}

\usepackage{hyperref}

\begin{document}
\title[The Capitulation Problem in Certain Pure Cubic Fields]
{The Capitulation Problem in Certain Pure Cubic Fields}

\author{Siham Aouissi}
\address{Number Theory and Information Security research team (TNSI)\\ Ecole Normale Sup\'erieure (ENS) of Moulay Ismail University (UMI)\\B.P. 3104\\Toulal\\Mekn\`es\\Morocco.}
\email{s.aouissi@umi.ac.ma}
\urladdr{https://sites.google.com/view/siham-aouissi}

\author{Daniel C. Mayer}
\address{Naglergasse 53\\8010 Graz\\Austria.}
\email{quantum.algebra@icloud.com}
\urladdr{http://www.algebra.at}

\thanks{Research of second author supported by the Austrian Science Fund (FWF): projects J0497-PHY, P26008-N25, and by the Research Executive Agency of the European Union (EUREA)}



\date{Sunday, 09 March 2025}

\begin{abstract}
Let \(\Gamma=\mathbb{Q}(\sqrt[3]{n})\) be a pure cubic field
with normal closure \(k=\mathbb{Q}(\sqrt[3]{n},\zeta)\),
where \(n>1\) denotes a cube free integer,
and \(\zeta\) is a primitive cube root of unity.
Suppose \(k\) possesses an
elementary bicyclic \(3\)-class group \(\mathrm{Cl}_3(k)\),
and the conductor of \(k/\mathbb{Q}(\zeta)\)
has the shape \(f\in\lbrace pq_1q_2,3pq,9pq\rbrace\)
where \(p\equiv 1\,(\mathrm{mod}\,9)\) and
\(q,q_1,q_2\equiv 2,5\,(\mathrm{mod}\,9)\) are primes.
It is disproved
that there are only two possible capitulation types \(\varkappa(k)\),
either type \(\mathrm{a}.1\), \((0000)\),
or type \(\mathrm{a}.2\), \((1000)\).
Evidence is provided, theoretically and experimentally,
of two further types,
\(\mathrm{b}.10\), \((0320)\), and
\(\mathrm{d}.23\), \((1320)\).

\end{abstract}

\maketitle

\hypersetup{linkcolor=blue}
\tableofcontents

\section{Introduction}
\label{s:Intro}

\noindent
Let \(n>1\) be a cube free integer, and
\(\Gamma=\mathbb{Q}(\sqrt[3]{n})\) be the pure cubic field with radicand \(n\),
i.e., a non-Galois field of degree three.
Denote by \(\zeta\) a primitive cube root of unity, and
let \(k_0=\mathbb{Q}(\zeta)\) be the third cyclotomic field.
Then the compositum \(k=\Gamma\cdot k_0\) is the normal closure of \(\Gamma\),
an absolutely dihedral field of degree six. 
Suppose \(k\) possesses an
elementary bicyclic \(3\)-class group
\(\mathrm{Cl}_3(k)\simeq (\mathbb{Z}/3\mathbb{Z})\times(\mathbb{Z}/3\mathbb{Z})\),
and the conductor of \(k/\mathbb{Q}(\zeta)\)
has the shape \(f\in\lbrace pq_1q_2,3pq,9pq\rbrace\)
where \(p\equiv 1\,(\mathrm{mod}\,9)\) and
\(q,q_1,q_2\equiv 2,5\,(\mathrm{mod}\,9)\) are primes.
In this situation, the \textit{capitulation type}
\(\varkappa(k)=(\ker(T_i))_{1\le i\le 4}\) of \(k\)
is defined by the kernels of the transfer homomorphisms
\(T_i:\mathrm{Cl}_3(k)\to\mathrm{Cl}_3(K_i)\) from \(k\)
to the four unramified cyclic cubic relative extensions
\(K_1,\ldots,K_4\) of \(k\).

In his Doctoral Thesis of 1992
\cite{Is1992},
Ismaili determined the possible capitulation types
of the normal closures of certain pure cubic fields.
In 
\cite{IsEM2005},
the possible capitulation types were narrowed down in 2005.
In the section \S\
\ref{s:Theoretical},
we show theoretically that some group theoretic arguments of \cite{IsEM2005},
used in the elimination process, were not justified,
and the purging went too far.
Our criticism is underpinned by experimental counter examples
in the sections \S\
\ref{s:Experimental}
and \S\
\ref{s:Coclass}.




{\tiny

\setlength{\unitlength}{1.0cm}
\begin{picture}(15,11)(-11,-10)

\put(-10,0.2){\makebox(0,0)[cb]{Degree}}

\put(-10,-2){\vector(0,1){2}}

\put(-10,-2){\line(0,-1){7}}
\multiput(-10.1,-2)(0,-1){8}{\line(1,0){0.2}}

\put(-10.2,-2){\makebox(0,0)[rc]{\(54\)}}
\put(-9.8,-2){\makebox(0,0)[lc]{Hilbert \(3\)-class field of normal closure and relative genus field}}
\put(-10.2,-4){\makebox(0,0)[rc]{\(18\)}}
\put(-9.8,-3.7){\makebox(0,0)[lb]{Unramified}}
\put(-9.8,-4){\makebox(0,0)[lc]{cyclic cubic}}
\put(-9.8,-4.3){\makebox(0,0)[lc]{extensions}}
\put(-10.2,-5){\makebox(0,0)[rc]{\(9\)}}
\put(-9.8,-5){\makebox(0,0)[lc]{Hilbert \(3\)-class fields}}
\put(-10.2,-6){\makebox(0,0)[rc]{\(6\)}}
\put(-9.8,-6){\makebox(0,0)[lc]{Dihedral normal closure}}
\put(-10.2,-7){\makebox(0,0)[rc]{\(3\)}}
\put(-9.8,-7){\makebox(0,0)[lc]{Conjugate pure cubic fields}}
\put(-10.2,-8){\makebox(0,0)[rc]{\(2\)}}
\put(-9.8,-8){\makebox(0,0)[lc]{Cyclotomic quadratic field}}
\put(-10.2,-9){\makebox(0,0)[rc]{\(1\)}}
\put(-9.8,-9){\makebox(0,0)[lc]{Rational base field}}

{\normalsize
\put(-3,0){\makebox(0,0)[cc]{Figure 1: Type III with \((k/k_0)^\ast=\mathrm{F}_3^1(k)\)}}
}



\put(0,-2){\circle*{0.2}}
\put(0.2,-2){\makebox(0,0)[lc]{\(\mathrm{F}_3^1(k)=(k/k_0)^\ast\)}}

\put(0,-2){\line(-3,-2){3}}
\put(0,-2){\line(-1,-2){1}}
\put(0,-2){\line(1,-2){1}}
\put(0,-2){\line(3,-2){3}}

\multiput(-3,-4)(2,0){4}{\circle*{0.2}}
\put(-3.3,-4){\makebox(0,0)[rc]{\(\mathrm{F}_3^1(\Gamma)\cdot k=\mathrm{F}_3^1(\Gamma^{\sigma})\cdot k=\mathrm{F}_3^1(\Gamma^{\sigma^2})\cdot k\)}}
\put(-2.7,-4){\makebox(0,0)[lc]{\(=K_4\)}}
\put(-1.2,-4){\makebox(0,0)[rc]{\(K_3\)}}
\put(1.2,-4){\makebox(0,0)[lc]{\(K_2\)}}
\put(3.2,-4){\makebox(0,0)[lc]{\(K_1\)}}

\put(0,-6){\line(-3,2){3}}
\put(0,-6){\line(-1,2){1}}
\put(0,-6){\line(1,2){1}}
\put(0,-6){\line(3,2){3}}

\put(0,-6){\circle*{0.2}}
\put(0.2,-6){\makebox(0,0)[lc]{\(k=\mathbb{Q}(\sqrt[3]{n},\zeta)\)}}

\put(0,-6){\line(0,-1){2}}
\put(0.2,-7){\makebox(0,0)[lc]{\(\mathrm{Gal}(k/k_0)=\langle\sigma\rangle\)}}

\put(0,-8){\circle*{0.2}}
\put(0.2,-8){\makebox(0,0)[lc]{\(k_0=\mathbb{Q}(\zeta)\)}}

\put(-5,-5){\line(2,1){2}}
\put(-2,-7){\line(2,1){2}}
\put(-2,-9){\line(2,1){2}}


\put(-5,-5){\circle{0.2}}
\put(-5.3,-5){\makebox(0,0)[rc]{\(\mathrm{F}_3^1(\Gamma)\),}}
\put(-4.7,-5){\makebox(0,0)[lc]{\(\mathrm{F}_3^1(\Gamma^{\sigma}),\mathrm{F}_3^1(\Gamma^{\sigma^2})\)}}

\put(-2,-7){\line(-3,2){3}}

\put(-2,-7){\circle{0.2}}
\put(-2.3,-7){\makebox(0,0)[rc]{\(\Gamma=\mathbb{Q}(\sqrt[3]{n})\),}}
\put(-1.7,-7){\makebox(0,0)[lc]{\(\Gamma^{\sigma},\Gamma^{\sigma^2}\)}}

\put(-2,-7){\line(0,-1){2}}

\put(-2,-9){\circle*{0.2}}
\put(-2.2,-9){\makebox(0,0)[rc]{\(\mathbb{Q}\)}}


\end{picture}

}

\section{Clarifications on capitulation types}
\label{s:Errors}

\noindent
The paper
\cite{IsEM2005}
states three theorems with erroneous proofs,
Th\'eor\`eme 2 in section 3 on page 53,
Th\'eor\`eme 3 in section 4 on page 56, and
Th\'eor\`eme 4 in section 5 on page 57.
Only for the first of these three theorems,
Th\'eor\`eme 2,
an explicit proof is established in
\cite[pp. 53--55]{IsEM2005},
which uses an incomplete version
of a group theoretic theorem,
Satz 14.17 on page 371,
by Huppert
\cite{Hu1967},
and consequently arrives
at the erroneous claims.
For Th\'eor\`eme 3,
 it is referred to
the proof of Th\'eor\`eme 2,
immediately after Th\'eor\`eme 3.
For Th\'eor\`eme 4,
it is also referred to
the proof of Th\'eor\`eme 2,
immediately before Proposition 8
and Th\'eor\`eme 4.
All three theorems concern
the capitulation kernels,
\(\varkappa(k)=(\ker(T_i:\mathrm{Cl}_3(k)\to\mathrm{Cl}_3(K_i)))_{1\le i\le 4}\),
of the \(3\)-class group
\(\mathrm{Cl}_3(k)\simeq (\mathbb{Z}/3\mathbb{Z})\times(\mathbb{Z}/3\mathbb{Z})\)
of the normal closure \(k=\mathbb{Q}(\sqrt[3]{n},\zeta)\)
of a pure cubic field \(\Gamma=\mathbb{Q}(\sqrt[3]{n})\)
in the four unramified cyclic cubic relative extensions
\(K_1,\ldots,K_4\) of \(k\).
The conductors of \(k/\mathbb{Q}(\zeta)\)
are \(f\in\lbrace pq_1q_2,3pq,9pq\rbrace\)
with primes \(p\equiv 1\,(\mathrm{mod}\,9)\) and
\(q,q_1,q_2\equiv 2,5\,(\mathrm{mod}\,9)\).
All three theorems
erroneously exclude two capitulation types
\(\varkappa(k)\in\lbrace (1320),(0320)\rbrace\),
based on the incomplete use of Huppert's theorem,
and state that only two capitulation types
\(\varkappa(k)\in\lbrace (1000),(0000)\rbrace\)
are possible.
In the particular situation where the auxiliary pure cubic field
\(\mathbb{Q}(\sqrt[3]{p})\) possesses a
differential principal factorization (DPF) of type \(\alpha\)
\cite[Thm. 2.1, p. 254]{AMITA2020},
the capitulation is even reduced to the single type
\(\varkappa(k)=(0000)\).


\section{Pure cubic fields of type III}
\label{s:TypeIII}

\noindent
In the last chapter of his Ph.D. Thesis
\cite[\S\ 3.4, pp. 44--59]{Is1992},
\cite{AAI2001},
Ismaili considers normal closures \(k\)
of pure cubic fields \(\Gamma\)
subject to the following conditions (see Figure 1):

\begin{enumerate}
\item
the class number \(h_\Gamma\) of \(\Gamma\) is exactly divisible by \(3\) and
the \(3\)-class group of \(k\) is elementary bicyclic,
\(\mathrm{Cl}_3(k)\simeq(\mathbb{Z}/3\mathbb{Z})\times(\mathbb{Z}/3\mathbb{Z})\),
\item
the relative genus field \((k/k_0)^\ast\)
coincides with the Hilbert \(3\)-class field \(\mathrm{F}_3^1(k)\) of \(k\),
\item
the conductor of \(k/k_0\)
is of the shape \(f\in\lbrace 3pq,9pq,pq_1q_2\rbrace\)
where \(p\equiv 1\,(\mathrm{mod}\,9)\) and
\(q,q_1,q_2\equiv 2,5\,(\mathrm{mod}\,9)\) are prime numbers.
\end{enumerate}
Ismaili calls normal closures \(k\) satisfying (1) and (2)
\textit{of type III}
\cite[Dfn. 3.1, p. 31]{Is1992},
and determines the class field theoretic \(3\)-extensions,
which are essential for the capitulation problem,
in the following way,
if (3) is satisfied additionally.
(We do not pay attention to the conductors \(f=9p\)
in Th\'eor\`eme 2, section 3, page 53 of
\cite{IsEM2005},
and \(f=pq\) with \(q\equiv 8\,(\mathrm{mod}\,9)\)
in Th\'eor\`eme 3, section 4, page 56, and we focus
on Th\'eor\`eme 4, section 5, page 57 of
\cite{IsEM2005}.)
If \(p\equiv 1\,(\mathrm{mod}\,3)\) is a prime number,
then \(p\) splits into a product \(p=\pi_1\pi_2\)
of two prime elements \(\pi_1\) and \(\pi_2\) in \(k_0\).
If such a prime \(p\) divides the radicand \(n\), then
\(p\) ramifies in \(k\) and there exist prime ideals
\(\mathfrak{P}_1,\mathfrak{P}_2\) of \(k\) such that
\begin{equation}
\label{eqn:Split}
\pi_1\mathcal{O}_k=\mathfrak{P}_1^3,\quad
\pi_2\mathcal{O}_k=\mathfrak{P}_2^3,\quad \text{ and } \quad
p\mathcal{O}_k=\mathfrak{P}_1^3\mathfrak{P}_2^3.
\end{equation}

\begin{theorem}
\label{thm:Ismaili}
The structure of all unramified abelian \(3\)-extensions of \(k\),
viewed as Kummer extensions, is given by adjoining
cube roots of the prime elements \(\pi_1,\pi_2\)
lying over the prime divisor \(p\equiv 1\,(\mathrm{mod}\,9)\)
of the conductor \(f\):
the relative genus field is bicyclic bicubic,
\begin{equation}
\label{eqn:Genus}
(k/k_0)^\ast=\mathrm{F}_3^1(k)=k(\sqrt[3]{\pi_1},\sqrt[3]{\pi_2}),
\end{equation}
and the four unramified cyclic cubic extensions of \(k\) are given by
\cite[Thm. 3.18, p. 58]{Is1992}
\begin{equation}
\label{eqn:Unramified}
K_1=k(\sqrt[3]{\pi_1\pi_2}),\quad
K_2=k(\sqrt[3]{\pi_2}),\quad
K_3=k(\sqrt[3]{\pi_1}),\quad \text{ and } \quad
K_4=k(\sqrt[3]{\pi_1\pi_2^2}).
\end{equation}
If the automorphism group \(\mathrm{Gal}(k/\Gamma)=\langle\tau\rangle\)
is generated by \(\tau\), then
\(\pi_2=\pi_1^\tau\),
the extensions \(K_2=K_3^\tau\) are conjugate non-Galois, 
and the extensions \(K_1\) and \(K_4\) are Galois over \(\mathbb{Q}\).
The associated norm class groups
\(H_i=\mathrm{Norm}_{K_i/k}\mathrm{Cl}_3(K_i)\), \(1\le i\le 4\),
are generated by the classes \(\lbrack\mathfrak{P}_j\rbrack\)
of prime ideals \(\mathfrak{P}_1,\mathfrak{P}_2=\mathfrak{P}_1^\tau\)
lying over the prime elements \(\pi_1,\pi_2\)
\cite[p. 51]{Is1992}:
\begin{equation}
\label{eqn:NormClass}
H_1=\langle\lbrack\mathfrak{P}_1\mathfrak{P}_2\rbrack\rangle,\quad
H_2=\langle\lbrack\mathfrak{P}_1\rbrack\rangle,\quad
H_3=\langle\lbrack\mathfrak{P}_2\rbrack\rangle,\quad \text{ and } \quad
H_4=\langle\lbrack\mathfrak{P}_1\mathfrak{P}_2^2\rbrack\rangle.
\end{equation}
\end{theorem}


\noindent
In order to present
the results of
\cite{Is1992}
on the
capitulation of \(3\)-classes of \(k\) in \(K_1,\ldots,K_4\),
i.e. the capitulation type \(\varkappa(k)\) of \(k\),
some background must be developed in Section
\ref{s:Capitulation}.


\section{Capitulation kernels and abelian type invariants}
\label{s:Capitulation}

\noindent
Since Scholz and Taussky
\cite[pp. 34--38]{SoTa1934}
only considered \textit{partial capitulation}
with cyclic kernels of order \(3\),
we use the capitulation types of Nebelung
\cite{Ne1989},
which are summarized in
\cite[Tbl. 6--7, pp. 492--493]{Ma2012}.
The latter also admit \textit{total capitulation}
with bicyclic kernel of order \(9\).
Let \(\mathcal{P}_K\) be the principal ideals and \(\mathcal{O}_K\) the integers 
of a number field \(K\). 


\begin{definition}
\label{dfn:TKT}
Let \(k\) be an algebraic number field
with elementary bicyclic \(3\)-class group
\(\mathrm{Cl}_3(k)\simeq(\mathbb{Z}/3\mathbb{Z})\times(\mathbb{Z}/3\mathbb{Z})\).
Denote by \(K_1,\ldots,K_4\) the four
unramified cyclic cubic extensions of \(k\),
and by \(T_i:\,\mathrm{Cl}_3(k)\to\mathrm{Cl}_3(K_i)\),
\(\mathfrak{a}\mathcal{P}_k\mapsto(\mathfrak{a}\mathcal{O}_{K_i})\mathcal{P}_{K_i}\),
the \textit{class extension homomorphisms} (briefly \textit{transfers})
from \(k\) to \(K_i\).
Then, the (unramified) \textit{capitulation type} (or \textit{transfer kernel type}, TKT) of \(k\)
\cite[Dfn. 1.1]{Ma2013}
is defined
in terms of the norm class groups of the extensions \(K_j/k\)
by \(\varkappa(k):=(\varkappa_1,\ldots,\varkappa_4)\) with
\begin{equation}
\label{eqn:TKT}
\varkappa_i:=
\begin{cases}
j & \text{ if } \ker(T_i)=\mathrm{Norm}_{K_j/k}\mathrm{Cl}_3(K_j), 
\text { where } 1\le i,j\le 4 \\
0 & \text{ if } \ker(T_i)=\mathrm{Cl}_3(k).
\end{cases}
\end{equation}
\end{definition}

\noindent
Here, \(\mathrm{Cl}_3(k)\) is the \(3\)-class group and \(\mathfrak{a}\) is an ideal of \(k\).
According to Hilbert's Theorem 94,
a transfer \(T_i\) cannot be injective,
that is, its kernel \(\ker(T_i)\) cannot be trivial. 


\begin{definition}
\label{dfn:PO}
Capitulation types
\(\varkappa=(\varkappa_1,\ldots,\varkappa_4)\) and \(\lambda=(\lambda_1,\ldots,\lambda_4)\)
are \textit{partially ordered} by the declaration
\(\varkappa\le\lambda\) \(:\Longleftrightarrow\)
\(\lbrack\varkappa_i\le\lambda_i\) component wise for all \(1\le i\le 4\rbrack\),
where, for each \(i=1,\ldots,4\), we put
\(\varkappa_i\le\lambda_i\) \(:\Longleftrightarrow\)
either \(\varkappa_i=\lambda_i\) or \(\lambda_i=0\).
\end{definition}

\noindent
The unique maximal capitulation type is a.1, \((0000)\),
but any capitulation type without zero components is minimal,
since it cannot shrink further, by Hilbert's Theorem 94.
In the sequel,
we are only concerned with four special capitulation types
\cite[pp. 492--493]{Ma2012},
subject to the following partial order.
Note that the types b.10,  \((0320)\), and a.2, \((1000)\), are incomparable.

\begin{equation}
\label{eqn:PO}
\begin{matrix}
                                 &          &          &  (0000), &  \mathrm{a}.1, \text{ maximal } \\
                                 &          & \nearrow & \uparrow &                                 \\
 \mathrm{b}.10,                  &  (0320)  &          &  (1000), &  \mathrm{a}.2                   \\
                                 & \uparrow & \nearrow &          &                                 \\
 \mathrm{d}.23, \text{ minimal } &  (1320)  &          &          &                                 \\
\end{matrix}
\end{equation}


\begin{definition}
\label{dfn:FP}
Let \(\varkappa=(\varkappa_1,\ldots,\varkappa_4)\) be a capitulation type.
A \textit{fixed point} of \(\varkappa\)
is a component \(\varkappa_i=i\), for some \(1\le i\le 4\).
A \textit{transposition} of \(\varkappa\)
is a pair of distinct components such that
\(\varkappa_i=j\) and \(\varkappa_j=i\),
for some \(1\le i\ne j\le 4\).
\end{definition}

\noindent
The types \(\mathrm{a}.2\), \((1000)\), and \(\mathrm{d}.23\), \((1320)\),
possess the fixed point \(1\).
The types \(\mathrm{b}.10\), \((0320)\), and \(\mathrm{d}.23\), \((1320)\),
contain a transposition \((2\mapsto 3\), \(3\mapsto 2)\).


Now we are in the position to state
the results
on the capitulation over normal closures \(k\)
of pure cubic fields \(\Gamma\) of type III
in
\cite[\S\ 3.4 (iii), Thm. 3.18, p. 58]{Is1992}.

\begin{theorem}
\label{thm:Capitulation}
If \(k\) is of type III,
there are four possibilities for the capitulation type,
\begin{equation}
\label{eqn:Capitulation}
\varkappa(k)\in\lbrace (1320),(0320),(1000),(0000)\rbrace,
\end{equation}
because, for exponents \(1\le a,b\le 2\), the power product of prime ideals
\begin{equation}
\label{eqn:Principal}
\mathfrak{P}_1^a\mathfrak{P}_2^b\mathcal{O}_{K_i}=\sqrt[3]{\pi_1^a\pi_2^b}\mathcal{O}_{K_i}
\text{ becomes principal for } i=
\begin{cases}
1 & \text{ if } a=b=1, \\
2 & \text{ if } a=0,\ b=1, \\
3 & \text{ if } a=1,\ b=0, \\
4 & \text{ if } a=1,\ b=2, \\
\end{cases}
\end{equation}
that is, the transfer kernels contain at least the classes
\begin{equation}
\label{eqn:Kernels}
\lbrack\mathfrak{P}_1\mathfrak{P}_2\rbrack\in\ker(T_1),\quad
\lbrack\mathfrak{P}_2\rbrack\in\ker(T_2),\quad
\lbrack\mathfrak{P}_1\rbrack\in\ker(T_3),\quad \text{ and } \quad
\lbrack\mathfrak{P}_1\mathfrak{P}_2^2\rbrack\in\ker(T_4),
\end{equation}
and since the prime ideal \(\mathfrak{p}\) lying over \(p\) in \(\Gamma\)
certainly becomes principal in the Hilbert \(3\)-class field \(\mathrm{F}_3^1(\Gamma)\),
the product \(\mathfrak{P}_1\mathfrak{P}_2=\mathfrak{p}\mathcal{O}_k\)
becomes principal in \(K_4=k\cdot\mathrm{F}_3^1(\Gamma)\),
whence \(\ker(T_4)=
\langle\lbrack\mathfrak{P}_1\mathfrak{P}_2^2\rbrack,
\lbrack\mathfrak{P}_1\mathfrak{P}_2\rbrack\rangle=\mathrm{Cl}_3(k)\).
Therefore, the minimal possible capitulation type is \(\varkappa(k)=(1320)\),
containing a fixed point and a transposition,
and since \(K_2\) and \(K_3\) are conjugate,
their transfer kernels may expand simultaneously only,
that is, \(\varkappa(k)=(1320)\) becomes \(\varkappa(k)=(1000)\).
However, if the fixed point expands,
then two other capitulation types, \(\varkappa(k)=(0320)\) and \(\varkappa(k)=(0000)\), arise.
\end{theorem}


\noindent
In the following corollary,
which describes the \textit{Galois structure} of
the relative \(3\)-genus field \((k/k_0)^\ast\),
we use the notation in angle brackets
\(\langle \mathrm{order},\mathrm{identifier}\rangle\)
of the SmallGroups library
\cite{BEO2005}
for finite groups of small order.

\begin{corollary}
\label{cor:Capitulation}
If \(k\) is of type III,
the absolute Galois group of the Hilbert \(3\)-class field of \(k\) is
\(G=\mathrm{Gal}(\mathrm{F}_3^1(k)/\mathbb{Q})\simeq\langle 54,13\rangle\).
Since the unramified cyclic cubic relative extensions
\(K_2\) and \(K_3\) of \(k\) are absolutely non-Galois conjugate,
their splitting field is the relative \(3\)-genus field \((k/k_0)^\ast=\mathrm{F}_3^1(k)\).
In the setting of Theorem
\ref{thm:Capitulation},
the other two unramified cyclic cubic relative extensions
\(K_1\) and \(K_4\) of \(k\) are absolutely Galois with non-abelian groups
\(\mathrm{Gal}(K_1/\mathbb{Q})\simeq\langle 18,4\rangle\) and
\(\mathrm{Gal}(K_4/\mathbb{Q})\simeq\langle 18,3\rangle\).
The maximal abelian normal subgroup of \(G\) is the relative \(3\)-genus group
\(\mathrm{Gal}((k/k_0)^\ast/k_0)\simeq\langle 27,5\rangle\)
which is elementary tricyclic.
\end{corollary}

\begin{proof}
Among the \(15\) groups of order \(54\),
\(G=\mathrm{Gal}((k/k_0)^\ast/\mathbb{Q})\simeq\langle 54,13\rangle\) is the unique group
with \(4\) triplets of conjugate subgroups \(\langle 18,3\rangle\) of order \(18\)
and a single self-conjugate subgroup \(\langle 18,4\rangle\) of order \(18\).
According to the Galois correspondence,
these subgroups are associated with
\(4\) triplets of conjugate pure cubic subfields,
and a single cyclic cubic subfield with conductor \(p\),
of the relative \(3\)-genus field \((k/k_0)^\ast\).
For the Dedekind species II, \(f=pq_1q_2\),
the pure cubic subfields are
one with conductor \(p\),
one with conductor \(q_1\cdot q_2\),
and two with conductor \(p\cdot q_1\cdot q_2\).
For the Dedekind species IB, \(f=3pq\), they are
one with conductor \(p\),
one with conductor \(3\cdot q\),
and two with conductor \(3\cdot p\cdot q\).
For the Dedekind species IA, \(f=9pq\), they are
one with conductor \(p\),
one with conductor \(3\cdot q\),
and two with conductor \(9\cdot p\cdot q\).
The group \(G=\langle 54,13\rangle\)
possesses a unique subgroup of order \(27\),
and this subgroup is abelian,
more precisely elementary tricyclic \(\langle 27,5\rangle\).
The extensions \(K_2\) and \(K_3\) are associated with
a pair of conjugate subgroups \(\langle 3,1\rangle\) of order \(3\),
and the extensions \(K_1\) and \(K_4\) are associated with
two self-conjugate subgroups \(\langle 3,1\rangle\) of order \(3\).
However, the quotients
\(\mathrm{Gal}(K_1/\mathbb{Q})\simeq\mathrm{Gal}((k/k_0)^\ast/\mathbb{Q})/\mathrm{Gal}((k/k_0)^\ast/K_1)\simeq\langle 18,4\rangle\) and
\(\mathrm{Gal}(K_4/\mathbb{Q})\simeq\mathrm{Gal}((k/k_0)^\ast/\mathbb{Q})/\mathrm{Gal}((k/k_0)^\ast/K_4)\simeq\langle 18,3\rangle\)
are different.
\end{proof}


\noindent
While
the Theorems
\ref{thm:Ismaili}
and
\ref{thm:Capitulation}
are perfectly correct,
the paper
\cite[Thm. 4, p. 57]{IsEM2005}
uses erroneous group theoretic arguments
to narrow down the four types \(\varkappa(k)\) to only two.
In fact, the paper \cite{IsEM2005} tries to eliminate the two capitulation types 
\(\mathrm{b}.10\), \((0320)\), and \(\mathrm{d}.23\), \((1320)\),
which contain a transposition \((2\mapsto 3\), \(3\mapsto 2)\).

This transposition can be characterized by
the structure of the \(3\)-class groups \(\mathrm{Cl}_3(K_i)\)
of the corresponding extensions with \(i\in\lbrace 2,3\rbrace\).


\begin{definition}
\label{dfn:ATI}
Let \(k\) be an algebraic number field
with elementary bicyclic \(3\)-class group
\(\mathrm{Cl}_3(k)\simeq(\mathbb{Z}/3\mathbb{Z})\times(\mathbb{Z}/3\mathbb{Z})\).
Denote by \(K_1,\ldots,K_4\) the four
unramified cyclic cubic extensions of \(k\).
Then, the \textit{abelian type invariants} (ATI) of the
four \(3\)-class groups \(\alpha(k):=(\mathrm{Cl}_3(K_i))_{1\le i\le 4}\)
are called the \textit{transfer target type} (TTT) of \(k\)
\cite[Dfn. 1.1]{Ma2013}.
Together with the transfer kernel type (TKT)
they form the \textit{Artin pattern} of \(k\)
\cite[Dfn. 4.3, p. 27]{Ma2015},
\begin{equation}
\label{eqn:Artin}
\mathrm{AP}(k):=\lbrack \alpha(k),\varkappa(k)\rbrack=
\lbrack (\mathrm{Cl}_3(K_i))_{i=1}^4, (\ker(T_i))_{i=1}^4\rbrack.
\end{equation}
Frequently, we shall write the ATI in logarithmic form,
e.g. \((111)\hat{=}(3,3,3)\), \((21)\hat{=}(9,3)\).
\end{definition}


In 2009, we have determined the
\textit{abelian quotient invariants} (AQI)
of the four maximal subgroups \(M_1,\ldots,M_4\),
that is,
the ATI of the commutator quotients \(M_i/M_i^\prime\), \(1\le i\le 4\),
of all finite \(3\)-groups \(G\) with elementary bicyclic \(G/G^\prime\)
\cite{Ma2014a}.
In particular, we proved the following remarkable theorem,
using absolute identifiers of the SmallGroups database
\cite{BEO2005}.

\begin{theorem}
\label{thm:Lo5Id3}
The group \(B:=\langle 243,3\rangle\) has
nilpotency class \(\mathrm{cl}(B)=3\),
coclass \(\mathrm{cc}(B)=2\),
soluble length \(\mathrm{sl}(B)=2\),
nuclear rank \(\nu(B)=2\),
relation rank \(d_2(B)=4\), and Artin pattern
\begin{equation}
\label{eqn:Lo5Id3}
\mathrm{AP}(B)=\lbrack (21,111,111,21),(0320)\rbrack.
\end{equation}
\begin{enumerate}
\item
Any finite \(3\)-group \(G\) with \(G/G^\prime\simeq (3,3)\)
and one of the \(\mathrm{TKT}\)s
\(\mathrm{b}.10\), \((0320)\),
\(\mathrm{d}.19\), \((2320)\),
\(\mathrm{d}.23\), \((1320)\),
or \(\mathrm{d}.25\), \((4320)\),
 is a descendant of \(B\).
\item
Any finite \(3\)-group \(G\) with \(G/G^\prime\simeq (3,3)\)
and coclass \(\mathrm{cc}(G)\ge 3\) is descendant of \(B\).
\item
Any descendant of \(B\) shares with \(B\) the stable part of the Artin pattern,
\(\lbrack\alpha_0,\varkappa_0\rbrack\)
with \(\alpha_0=(\ldots,111,111,\ldots)\) and \(\varkappa_0=(\ldots 32\ldots)\),
in particular, the conspicuous transposition \((2\mapsto 3\), \(3\mapsto 2)\).
\end{enumerate}
\end{theorem}

\begin{proof}
The bifurcation group \(B\)
with nuclear rank \(\nu(B)=2\)
is determined uniquely by its identifier
in the SmallGroups database
\cite{BEO2005}.
Its AQI were determined in
\cite[\S\ 4.3, Tbl. 4.3, p. 434]{Ma2014a}.
\(B\) is the smallest group with TKT \(\mathrm{b}.10\), \((0320)\).
The types
\(\mathrm{d}.19\), \((2320)\),
\(\mathrm{d}.23\), \((1320)\), and
\(\mathrm{d}.25\), \((4320)\),
arise from type \(\mathrm{b}.10\), \((0320)\),
by shrinking of the \textit{polarization} (first component).
The AQI of the infinitely many descendants of \(B\)
were determined in
\cite[\S\ 4.4, Tbl. 4.5, p. 438]{Ma2014a}
and
\cite[\S\ 4.5, Tbl. 4.7, p. 441]{Ma2014a}
for coclass \(\mathrm{cc}(G)=2\), and in
\cite[\S\ 4.6, Thm. 4.5, pp. 444--445]{Ma2014a}
for the bifurcation to coclass \(\mathrm{cc}(G)\ge 3\).
The statement concerning coclass \(\mathrm{cc}(G)\ge 3\)
is due to Nebelung
\cite{Ne1989},
who discovered the crucial bifurcation from coclass \(2\) to coclass \(3\).
Since the components \((32)\) of the TKT cannot shrink,
according to Hilbert's Theorem 94,
they remain stable for any descendant of \(\langle 243,3\rangle\).
Accordingly, the corresponding components \((111,111)\) of the TTT (AQI)
cannot expand and form the \textit{stabilization} of the Artin pattern.
\end{proof}


\section{Little and big two-stage towers}
\label{s:LittleAndBig}

\noindent
\begin{definition}
\label{dfn:LittleAndBig}
Let \(k\) be an algebraic number field
with elementary bicyclic \(3\)-class group
\(\mathrm{Cl}_3(k)\simeq(\mathbb{Z}/3\mathbb{Z})\times(\mathbb{Z}/3\mathbb{Z})\).
Denote by \(K_1,\ldots,K_4\) the four
unramified cyclic cubic extensions of \(k\).
Then, the automorphism groups
\(\mathfrak{G}_i:=\mathrm{Gal}(\mathrm{F}_3^1(K_i)/k)\), \(1\le i\le 4\),
of the first Hilbert \(3\)-class fields of \(K_i\)
are called the groups of the four \textit{little two-stage towers} of \(k\),
and the automorphism group
\(G:=\mathrm{Gal}(\mathrm{F}_3^2(k)/k)\)
of the second Hilbert \(3\)-class field of \(k\)
is called the group of the \textit{big two-stage tower} of \(k\).
\end{definition}


\begin{theorem}
\label{thm:LittleAndBig}
Let the assumptions in Definition
\ref{dfn:LittleAndBig}
be satisfied.
Generally, the connection between groups of little and big two-stage towers is given by
\begin{equation}
\label{eqn:LittleAndBig}
\mathfrak{G}_i\simeq G/M_i^\prime, \qquad \mathfrak{A}_i\simeq\mathrm{Cl}_3(K_i)\simeq M_i/M_i^\prime
\end{equation}
where \(M_1,\ldots,M_4\) denote the four maximal subgroups of \(G\),
and \(\mathfrak{A}_i:=\mathrm{Gal}(\mathrm{F}_3^1(K_i)/K_i)\) is the distinguished
abelian maximal subgroup of \(\mathfrak{G}_i\), for \(1\le i\le 4\).
The commutator quotients
\(G/G^\prime\simeq\mathfrak{G}_i/\mathfrak{G}_i^\prime\simeq\mathrm{Cl}_3(k)
\simeq(\mathbb{Z}/3\mathbb{Z})\times(\mathbb{Z}/3\mathbb{Z})\)
are isomorphic, here elementary bicyclic.

In particular, if \(K_i\) has \(3\)-class group
\(\mathrm{Cl}_3(K_i)\simeq
(\mathbb{Z}/3\mathbb{Z})\times(\mathbb{Z}/3\mathbb{Z})\times(\mathbb{Z}/3\mathbb{Z})\),
the elementary tricyclic \(3\)-group,
for some \(i=1,\ldots,4\),
then the little two-stage tower has the group
with SmallGroups identifier
\cite{BEO2005}
\begin{equation}
\label{eqn:Syl3A9}
\mathfrak{G}_i\simeq\langle 81,7\rangle\simeq\mathrm{Syl}_3(A_9).
\end {equation}
\end{theorem}

\begin{proof}
For the proof of the general statements, we refer to
\cite[\S\ 2, Prop. 2.1, pp. 417--419]{Ma2014a}.
When the unramified cyclic cubic extension \(K_i/k\)
possesses an elementary tricyclic \(3\)-class group
\(\mathrm{Cl}_3(K_i)\simeq (3,3,3)\), then 
\(\lbrack\mathrm{F}_3^1(K_i):K_i\rbrack=3^3=27\)
is the relative degree of its Hilbert \(3\)-class field, and
the group \(\mathfrak{G}_i\) of the little two-stage tower
is of maximal class and order
\(\#\mathfrak{G}_i=\lbrack\mathrm{F}_3^1(K_i):K_i\rbrack\cdot\lbrack K_i:k\rbrack=3^3\cdot 3=81\),
according to
\cite[Satz 4, Korollar, p. 9]{HeSm1982},
because \(\mathrm{Cl}_3(k)\simeq (3,3)\).
Since \textit{the distinguished abelian maximal subgroup}
\(\mathfrak{A}_i=\mathrm{Gal}(\mathrm{F}_3^1(K_i)/K_i)\simeq\mathrm{Cl}_3(K_i)\)
of \(\mathfrak{G}_i\) is \textit{elementary tricyclic},
\(\mathfrak{G}_i\) is uniquely determined  as
\(\langle 81,7\rangle\simeq\mathrm{Syl}_3(A_9)\),
the Sylow \(3\)-subgroup of the alternating group \(A_9\) of degree \(9\),
according to
\cite[\S\ 4.1, Thm. 4.1, p. 427]{Ma2014a}.
\end{proof}


\section{Theoretical proof }
\label{s:Theoretical}

\noindent
We state the incorrect hypothesis
\cite[Thm. 4, p. 57]{IsEM2005}
in a form which is more detailed than in the original paper.
However, we exclude that \(q_i\equiv 8\,(\mathrm{mod}\,9)\),
for \(i=1\) or \(i=2\).

\begin{hypothesis}
\label{hyp:ElMesaoudi}
Let \(k=\mathbb{Q}(\sqrt[3]{n},\zeta)\) be the normal closure
of a pure cubic number field \(\Gamma=\mathbb{Q}(\sqrt[3]{n})\).
Suppose the cubefree radicand \(n\)
is a positive integer with one of the following decompositions
in prime numbers \(p\equiv 1\,(\mathrm{mod}\,9)\) and \(q,q_1,q_2\equiv 2,5\,(\mathrm{mod}\,9)\):
\begin{equation}
\label{eqn:Radicands}
\begin{aligned}
n &= p^eq_1^{e_1}q_2^{e_2}\equiv\pm 1\,(\mathrm{mod}\,9)
\text{ with integers } e,e_1,e_2\in\lbrace 1,2\rbrace \text{ (species II)}, \\
n &= p^{e_1}q^{e_2}\not\equiv\pm 1\,(\mathrm{mod}\,9)
\text{ with integers } e_1,e_2\in\lbrace 1,2\rbrace \text{ (species IB)}, \\
n &= 3^ep^{e_1}q^{e_2}
\text{ with integers } e,e_1,e_2\in\lbrace 1,2\rbrace \text{ (species IA)}.
\end{aligned}
\end{equation}
Assume the \(3\)-class group of \(k\) is elementary bicyclic,
\(\mathrm{Cl}_3(k)\simeq(\mathbb{Z}/3\mathbb{Z})\times(\mathbb{Z}/3\mathbb{Z})\),
the \(3\)-class group of \(\Gamma\) is elementary cyclic,
\(\mathrm{Cl}_3(\Gamma)\simeq\mathbb{Z}/3\mathbb{Z}\),
and the relative \(3\)-genus field \((k/k_0)^\ast\)
coincides with 
the Hilbert \(3\)-class field \(\mathrm{F}_3^1(k)\),
that is, \(k\) is of
type III, as drawn in Figure \(1\).
Then there are \textbf{only two} possibilities for the capitulation type of \(k\),
\begin{equation}
\label{eqn:Reduction}
\varkappa(k)\in\lbrace (1000),(0000)\rbrace,
\end{equation}
in the four unramified cyclic cubic relative extensions
\(K_1,\ldots,K_4\) of \(k\).
If the index of the subgroup \(U_0\) generated by all subfield units
in the group \(U(\Gamma_p)\)
of the pure cubic auxiliary field \(\Gamma_p=\mathbb{Q}(\sqrt[3]{p})\)
is equal to \((U(\Gamma_p):U_0)=1\),
that is, \(\Gamma_p\) is of differential principal factorization type \(\alpha\),
then \textbf{only the single} capitulation type \(\varkappa(k)=(0000)\) is possible.
\end{hypothesis}

\begin{warning}
\label{wrn:ElMesaoudi}
The proof of 
\cite[\S\ 5, Thm. 4, p. 57]{IsEM2005}
was conducted simultaneously with the proofs of 
\cite[\S\ 4, Thm. 3, p. 56]{IsEM2005}
and
\cite[\S\ 3, Thm. 2, p. 53]{IsEM2005}
in
\cite[\S\ 3, pp. 53--55]{IsEM2005},
making use of a theorem of Huppert
which is cited incompletely as 
\cite[\S\ 2, Prop. 4, p. 51]{IsEM2005}.
Although we only found explicit numerical counter examples to
\cite[\S\ 5, Thm. 4, p. 57]{IsEM2005}
in Section
\ref{s:Experimental},
it must be pointed out that
none of the above mentioned Theorems 2,3,4 in
\cite{IsEM2005}
was proved reliably,
and violations of Theorems 2,3 could still exist.
\end{warning}


\noindent
Now we explain how the erroneous Hypothesis
\ref{hyp:ElMesaoudi}
can be disproved theoretically.
For this purpose,
we first state the correct form
of the group theoretic result
 that was used in \cite{IsEM2005}
to derive the incorrect Hypothesis
\ref{hyp:ElMesaoudi}.

\noindent
Let \(p\) be a prime number and
\(G\) be a finite \(p\)-group with
lower central series \(\gamma_2{G}=\lbrack G,G\rbrack\) and
\(\gamma_i{G}=\lbrack\gamma_{i-1}{G},G\rbrack\) for \(i\ge 3\).
Denote by \(\gamma_1{G}\) the \textit{two-step centralizer} of \(G\) such that
\(\gamma_1{G}/\gamma_4{G}\) is the centralizer of
\(\gamma_2{G}/\gamma_4{G}\) in \(G/\gamma_4{G}\).
The correct form of 
\cite[\S\ 2, Prop. 4, p. 51]{IsEM2005}
is the following theorem by Huppert
\cite[Satz 14.17, p. 371]{Hu1967}.

\begin{theorem}
\label{thm:Huppert}
If \(G\) is a \(3\)-group of maximal nilpotency class
\textbf{with order} \(\mathbf{\#G=3^n}\), \(\mathbf{n>4}\), \textbf{that is}, \(\mathbf{\#G\ge 243}\),
then \(G\) is metabelian, i.e., \(G^{\prime\prime}=1\), and
\(\gamma_1{G}\) is a metacyclic group with nilpotency class \(\mathrm{cl}(\gamma_1{G})\le 2\).
\end{theorem}

\begin{proof}
This theorem is
\cite[Satz 14.17, p. 371]{Hu1967}.
Its proof makes use of 
\cite[Satz 14.16, p. 370]{Hu1967},
where statements about a \(p\)-group \(G\)
of maximal nilpotency class
\textbf{with order} \(\mathbf{\#G=p^n}\), \(\mathbf{n>p+1}\),
are proved,
\textbf{that is}, \(\mathbf{n>3+1=4}\), \textbf{if} \(\mathbf{p=3}\).
\end{proof}


\noindent
Now we state and prove the correct theorem
which must replace the erroneous Hypothesis
\ref{hyp:ElMesaoudi}.
We prefer to express the ramification
in terms of the conductor \(f\)
rather than the radicand \(n\)
with cumbersome and useless exponents \(e,e_1,e_2\).

\begin{theorem}
\label{thm:Correction}
Let \(k=\mathbb{Q}(\sqrt[3]{n},\zeta)\) be the normal closure
of a pure cubic number field \(\Gamma=\mathbb{Q}(\sqrt[3]{n})\).
Suppose the conductor \(f\) of \(k/k_0\)
is a positive integer with one of the following decompositions
in prime numbers \(p\equiv 1\,(\mathrm{mod}\,9)\) and \(q,q_1,q_2\equiv 2,5\,(\mathrm{mod}\,9)\):
\begin{equation}
\label{eqn:Conductors}
\begin{aligned}
f &= pq_1q_2
\text{ (species II)}, \\
f &= 3pq
\text{ (species IB)}, \\
f &= 9pq
\text{ (species IA)}.
\end{aligned}
\end{equation}
Assume the \(3\)-class group of \(k\) is elementary bicyclic,
\(\mathrm{Cl}_3(k)\simeq(\mathbb{Z}/3\mathbb{Z})\times(\mathbb{Z}/3\mathbb{Z})\),
the \(3\)-class group of \(\Gamma\) is elementary cyclic,
\(\mathrm{Cl}_3(\Gamma)\simeq\mathbb{Z}/3\mathbb{Z}\),
and the relative \(3\)-genus field \((k/k_0)^\ast\)
coincides with 
the Hilbert \(3\)-class field \(\mathrm{F}_3^1(k)\),
that is, \(k\) is of Ismaili's type III, as drawn in Figure \(1\).
Then, up to equivalence,
there are \textbf{four} possibilities for the capitulation type of \(k\)
in the four unramified cyclic cubic relative extensions
\(K_1,\ldots,K_4\) of \(k\),
either
\begin{equation}
\label{eqn:Total}
\varkappa(k)\in\lbrace (1000),(0000)\rbrace,
\end{equation}
if the unit norm index is \((U(k):\mathrm{Norm}_{K_i/k}U(K_i))=3\), for \(i\in\lbrace 2,3\rbrace\),
or
\begin{equation}
\label{eqn:Partial}
\varkappa(k)\in\lbrace (1320),(0320)\rbrace,
\end{equation}
if the unit norm index is \((U(k):\mathrm{Norm}_{K_i/k}U(K_i))=1\), for \(i\in\lbrace 2,3\rbrace\).

If the index of the subgroup \(U_0\) generated by all subfield units
in the group \(U(\Gamma_p)\)
of the pure cubic auxiliary field \(\Gamma_p=\mathbb{Q}(\sqrt[3]{p})\)
is equal to \((U(\Gamma_p):U_0)=1\),
that is, \(\Gamma_p\) is of DPF-type \(\alpha\),
then \textbf{only two} capitulation types are possible, depending on the unit norm index, \\
either \(\varkappa(k)=(0000)\), if  \((U(k):\mathrm{Norm}_{K_i/k}U(K_i))=3\), for \(i=2,3\), \\
or \(\varkappa(k)=(0320)\), if \((U(k):\mathrm{Norm}_{K_i/k}U(K_i))=1\), for \(i=2,3\).
\end{theorem}

\begin{proof} (of Theorem
\ref{thm:Correction} by disproof of Hypothesis
\ref{hyp:ElMesaoudi}.)
In \cite{IsEM2005},
 it is attempted to discourage the capitulation types
\(\mathrm{b}.10\), \((0320)\), and
\(\mathrm{d}.23\), \((1320)\),
by proving that in Formula
\eqref{eqn:Kernels}
not only
\(\lbrack\mathfrak{P}_2\rbrack\in\ker(T_2)\)
but also
\(\lbrack\mathfrak{P}_1\rbrack\in\ker(T_2)\)
\cite[\S\ 3, Proof of Thm. 2, item (1), p. 53]{IsEM2005},
and thus
\(\ker(T_3)\simeq\ker(T_2)=
\langle\lbrack\mathfrak{P}_1\rbrack,\lbrack\mathfrak{P}_2\rbrack\rangle=\mathrm{Cl}_3(k)\).
They use contraposition,
trying to derive a contradiction from the assumption that
\(\lbrack\mathfrak{P}_1\rbrack\not\in\ker(T_2)\)
\cite[\S\ 3, Proof of Thm. 2, item (1), p. 54]{IsEM2005},
which implies 
\(\mathrm{Cl}_3(K_3)\simeq\mathrm{Cl}_3(K_2)\simeq (3,3,3)\),
according to Theorem
\ref{thm:Lo5Id3}
and Formula
\eqref{eqn:Lo5Id3},
and consequently
\(\mathfrak{G}_2\simeq\mathfrak{G}_3\simeq\langle 81,7\rangle\simeq\mathrm{Syl}_3(A_9)\)
with \(\#\mathfrak{G}_2=\#\mathfrak{G}_3=81\),
according to Theorem
\ref{thm:LittleAndBig},
Formula
\eqref{eqn:Syl3A9}.

Then,
\textbf{the incomplete theorem of Huppert
\cite[\S\ 2, Prop. 4, p. 51]{IsEM2005}  is applied}
without the crucial condition \(\#G\ge 243\),
to the groups \(\mathfrak{G}_2\) and \(\mathfrak{G}_3\) of the little two-stage towers
\cite[\S\ 3, Proof of Thm. 2, item (1)(iv), p. 55]{IsEM2005},
and conclude that the two-step centralizers
\(\mathfrak{A}_2=\gamma_1\mathfrak{G}_2\) and \(\mathfrak{A}_3=\gamma_1\mathfrak{G}_3\) must be metacyclic,
in contradiction to \(\mathfrak{A}_2\simeq \mathfrak{A}_3\simeq (3,3,3)\).
However, the correct theorem of Huppert
cannot be applied to a group of order \(3^4=81\), 
and thus \textbf{no contradiction arises}
when the capitulation kernel \(\ker(T_2)\) is cyclic of order \(3\),
generated by \(\lbrack\mathfrak{P}_2\rbrack\),
i.e., \(\varkappa(k)=(\varkappa_1 320)\) with \(\varkappa_1\in\lbrace 0,1\rbrace\).

Finally, the well-known Galois-cohomological theorem on the Herbrand quotient
of the \textit{unramified} cyclic relative extension \(K_i/k\) of \textit{odd} prime degree \(\lbrack K_i:k\rbrack=3\),
\(\#\ker(T_i)=\#\mathrm{H}^1(\langle\sigma_i\rangle,U(K_i))
=\lbrack K_i:k\rbrack\cdot\#\hat{\mathrm{H}}^0(\langle\sigma_i\rangle,U(K_i))
=3\cdot (U(k):\mathrm{Norm}_{K_i/k}U(K_i))\), for \(1\le i\le 4\),
where \(\langle\sigma_i\rangle=\mathrm{Gal}(K_i/k)\simeq\mathbb{Z}/3\mathbb{Z}\),
admits the distinction between total and partial capitulation in \(K_2\) and \(K_3\): \\
\(\varkappa(k)\in\lbrace (1000),(0000)\rbrace\) if and only if \((U(k):\mathrm{Norm}_{K_i/k}U(K_i))=3\), for \(i=2,3\), and \\
\(\varkappa(k)\in\lbrace (1320),(0320)\rbrace\) if and only if \((U(k):\mathrm{Norm}_{K_i/k}U(K_i))=1\), for \(i=2,3\).
\end{proof}


\section{Experimental proof }
\label{s:Experimental}

\noindent
In this concluding section,
we underpin our theoretical disproof of
the erroneous claims in
\cite[Thm. 4, p. 57]{IsEM2005}
by extensive numerical counter examples
under the GRH,
compiled with the aid of class field theoretic routines by Fieker
\cite{Fi2001},
which are implemented in the computer algebra system Magma
\cite{BCP1997}, \cite{BCFS2023}, \cite{MAGMA2024}.
Since the normal closures \(k=\mathbb{Q}(\sqrt[3]{n},\zeta)\)
of pure cubic fields \(\Gamma=\mathbb{Q}(\sqrt[3]{n})\)
are constructed as \(3\)-ring class fields modulo \(3\)-admissible conductors \(f\),
they are automatically produced as \textit{multiplets} \((k_1,\ldots,k_m)\)
of pairwise non-isomorphic fields sharing the common conductor \(f\)
with multiplicities \(m=m(f)\) determined by the formulas in
\cite[Thm. 2.1, p. 833]{Ma1992}
or
\cite[Cor. 3.2, p. 2219]{Ma2014b}.
In Table
\ref{tbl:PQQ122}--\ref{tbl:PQQ125},
we show that Th\'eor\`eme 4 is violated by
conductors \(f=pq_1q_2\) of Dedekind species II
with primes \(p\equiv 1\,(\mathrm{mod}\,9)\) 
and various combinations of
\(q_1,q_2\equiv 2,5\,(\mathrm{mod}\,9)\).
Columns begin with prime factors \(p,q_1,q_2\) and conductor \(f\),
continue with normalized radicand \(n_i\), Artin pattern consisting of
capitulation type (CT) \(\varkappa_i\) and logarithmic abelian type invariants \(\alpha_i\),
and a minimal differential principal factor (DPF)
(norm of an ambiguous principal ideal),
split into \(m\) rows, according to the multiplicity \(m=m(f)\) of the conductor,
and end with the (DPF-)type of the pure cubic auxiliary field \(\mathbb{Q}(\sqrt[3]{p})\),
which decides
whether the \(\varkappa_i\) may be of CT \(\mathrm{d}.23\) (for type \(\gamma\))
or not (for type \(\alpha\)).
The radicand \(n_i=a\cdot b^2\) is called \textit{normalized},
if \(1\le b<a\) are square free coprime integers.


\subsection{Counterexamples of Dedekind species II}
\label{ss:DedekindII}

First, we list lots of counterexamples
of Dedekind species II
with conductors of the shape \(3\nmid f=pq_1q_2\),
giving rise to \textit{doublets},
in Table
\ref{tbl:PQQ122},
\ref{tbl:PQQ155}, and
\ref{tbl:PQQ125}.

Only in the mixed case
\(q_1\equiv 2\,(\mathrm{mod}\,9)\), \(q_2\equiv 5\,(\mathrm{mod}\,9)\)
of Table
\ref{tbl:PQQ125}
with most numerous hits,
there occurs a single instance of capitulation type (CT) d.23, \(\varkappa\sim (1320)\),
because the \textit{pure cubic auxiliary field} \(\mathbb{Q}(\sqrt[3]{p})\) is of DPF type \(\gamma\).
The primes \(p\equiv 1\,(\mathrm{mod}\,9)\) with type \(\gamma\)
are collected in the sequence A363717
of the On-line Encyclopedia of Integer Sequences (OEIS)
\cite{OEIS2025}.
They start with \(541,919\).


\renewcommand{\arraystretch}{1.1}

\begin{table}[ht]
\caption{Doublets, \(m=2\), with \(f=pq_1q_2\), \(q_1\equiv 2\,(\mathrm{mod}\,9)\), \(q_2\equiv 2\,(\mathrm{mod}\,9)\)}
\label{tbl:PQQ122}
\begin{center}

{\small

\begin{tabular}{|rrr|r|ccccc|c|}
\hline
   \(p\) & \(q_1\) & \(q_2\) &    \(f\) &                          \(n_i\) &   CT & \(\varkappa_i\) &      \(\alpha_i\) &                DPF &       Type \\
\hline
  \(19\) &   \(2\) &  \(11\) &  \(418\) &      \(836=2^2\cdot 11\cdot 19\) & b.10 &        \(4001\) & \(111,21,21,111\) &             \(11\) & \(\alpha\) \\
         &         &         &          &     \(4598=2\cdot 11^2\cdot 19\) & b.10 &        \(4001\) & \(111,21,21,111\) &             \(11\) &            \\
  \(19\) &   \(2\) &  \(83\) & \(3154\) & \(119852=2^2\cdot 19^2\cdot 83\) & b.10 &        \(3010\) & \(111,21,111,21\) &   \(38=2\cdot 19\) & \(\alpha\) \\
         &         &         &          &     \(6308=2^2\cdot 19\cdot 83\) & b.10 &        \(0043\) & \(21,21,111,111\) & \(76=2^2\cdot 19\) &            \\
  \(37\) &   \(2\) &  \(11\) &  \(814\) &     \(8954=2\cdot 11^2\cdot 37\) & b.10 &        \(3010\) & \(111,21,111,21\) &             \(11\) & \(\alpha\) \\
         &         &         &          &     \(1628=2^2\cdot 11\cdot 37\) & b.10 &        \(0043\) & \(21,21,111,111\) &             \(11\) &            \\
  \(37\) &   \(2\) &  \(29\) & \(2146\) &     \(4292=2^2\cdot 29\cdot 37\) & b.10 &        \(0320\) & \(21,111,111,21\) &             \(29\) & \(\alpha\) \\
         &         &         &          &    \(62234=2\cdot 29^2\cdot 37\) & b.10 &        \(3010\) & \(111,21,111,21\) &             \(29\) &            \\
  \(37\) &   \(2\) &  \(47\) & \(3478\) &     \(6956=2^2\cdot 37\cdot 47\) & b.10 &        \(3010\) & \(111,21,111,21\) &             \(47\) & \(\alpha\) \\
         &         &         &          &   \(163466=2\cdot 37\cdot 47^2\) & b.10 &        \(3010\) & \(111,21,111,21\) &             \(47\) &            \\
 \(109\) &   \(2\) &  \(11\) & \(2398\) &   \(26378=2\cdot 11^2\cdot 109\) & b.10 &        \(0402\) & \(21,111,21,111\) &              \(2\) & \(\alpha\) \\
         &         &         &          &    \(4796=2^2\cdot 11\cdot 109\) & b.10 &        \(0320\) & \(21,111,111,21\) &              \(2\) &            \\
 \(127\) &   \(2\) &  \(11\) & \(2794\) &    \(5588=2^2\cdot 11\cdot 127\) & b.10 &        \(0320\) & \(21,111,111,21\) &              \(2\) & \(\alpha\) \\
         &         &         &          &   \(30734=2\cdot 11^2\cdot 127\) & b.10 &        \(3010\) & \(111,21,111,21\) &              \(2\) &            \\
 \(199\) &   \(2\) &  \(11\) & \(4378\) &    \(8756=2^2\cdot 11\cdot 199\) & b.10 &        \(0402\) & \(32,111,22,111\) &             \(11\) & \(\alpha\) \\
         &         &         &          &   \(48158=2\cdot 11^2\cdot 199\) & b.10 &        \(0320\) & \(32,111,111,22\) &             \(11\) &            \\
\hline
\end{tabular}

}

\end{center}
\end{table}


\renewcommand{\arraystretch}{1.1}

\begin{table}[ht]
\caption{Doublets, \(m=2\), with \(f=pq_1q_2\), \(q_1\equiv 5\,(\mathrm{mod}\,9)\), \(q_2\equiv 5\,(\mathrm{mod}\,9)\)}
\label{tbl:PQQ155}
\begin{center}

{\small

\begin{tabular}{|rrr|r|ccccc|c|}
\hline
   \(p\) & \(q_1\) & \(q_2\) &    \(f\) &                         \(n_i\) &   CT & \(\varkappa_i\) &      \(\alpha_i\) &                DPF &        Type \\
\hline
  \(37\) &  \(23\) &   \(5\) & \(4255\) & \(97865=5^2\cdot 23^2\cdot 37\) & b.10 &        \(3010\) & \(111,21,111,21\) &             \(23\) & \(\alpha\) \\
         &         &         &          &   \(21275=5^2\cdot 23\cdot 37\) & b.10 &        \(3010\) & \(111,21,111,21\) &             \(23\) &            \\
\hline
\end{tabular}

}

\end{center}
\end{table}


\renewcommand{\arraystretch}{1.1}

\begin{table}[ht]
\caption{Doublets, \(m=2\), with \(f=pq_1q_2\), \(q_1\equiv 2\,(\mathrm{mod}\,9)\), \(q_2\equiv 5\,(\mathrm{mod}\,9)\)}
\label{tbl:PQQ125}
\begin{center}

{\small

\begin{tabular}{|rrr|r|ccccc|c|}
\hline
   \(p\) & \(q_1\) & \(q_2\) &    \(f\) &                         \(n_i\) &   CT & \(\varkappa_i\) &      \(\alpha_i\) &                DPF &       Type \\
\hline
  \(19\) &   \(2\) & \(113\) & \(4294\) &  \(81586=2\cdot 19^2\cdot 113\) & b.10 &        \(0043\) & \(21,21,111,111\) & \(76=2^2\cdot 19\) & \(\alpha\) \\
         &         &         &          &     \(4294=2\cdot 19\cdot 113\) & b.10 &        \(3010\) & \(111,21,111,21\) &   \(38=2\cdot 19\) &            \\
  \(19\) &  \(11\) &   \(5\) & \(1045\) &      \(1045=5\cdot 11\cdot 19\) & b.10 &        \(3010\) & \(111,21,111,21\) &             \(11\) & \(\alpha\) \\
         &         &         &          &   \(19855=5\cdot 11\cdot 19^2\) & b.10 &        \(4001\) & \(111,21,21,111\) &             \(11\) &            \\
  \(19\) &  \(11\) &  \(23\) & \(4807\) &     \(4807=11\cdot 19\cdot 23\) & b.10 &        \(4001\) & \(111,21,21,111\) &             \(11\) & \(\alpha\) \\
         &         &         &          &  \(91333=11\cdot 19^2\cdot 23\) & b.10 &        \(3010\) & \(111,21,111,21\) &             \(11\) &            \\
  \(37\) &   \(2\) &  \(23\) & \(1702\) &      \(1702=2\cdot 23\cdot 37\) & b.10 &        \(3010\) & \(111,21,111,21\) &             \(23\) & \(\alpha\) \\
         &         &         &          &   \(62974=2\cdot 23\cdot 37^2\) & b.10 &        \(4001\) & \(111,21,21,111\) &             \(23\) &            \\
  \(37\) &  \(11\) &   \(5\) & \(2035\) &      \(2035=5\cdot 11\cdot 37\) & b.10 &        \(4001\) & \(111,21,21,111\) &             \(11\) & \(\alpha\) \\
         &         &         &          &   \(75295=5\cdot 11\cdot 37^2\) & b.10 &        \(3010\) & \(111,21,111,21\) &             \(11\) &            \\
  \(37\) &   \(29\) &  \(5\) & \(5365\) &      \(5365=5\cdot 29\cdot 37\) & b.10 &        \(0043\) & \(21,21,111,111\) &             \(29\) & \(\alpha\) \\
         &         &         &          &  \(198505=5\cdot 29\cdot 37^2\) & b.10 &        \(3010\) & \(111,21,111,21\) &             \(29\) &            \\
 \(109\) &   \(2\) &   \(5\) & \(1090\) & \(10900=2^2\cdot 5^2\cdot 109\) & b.10 &        \(0320\) & \(21,111,111,21\) &              \(2\) & \(\alpha\) \\
         &         &         &          &      \(1090=2\cdot 5\cdot 109\) & b.10 &        \(0043\) & \(21,21,111,111\) &              \(2\) &            \\
 \(163\) &   \(2\) &   \(5\) & \(1630\) &      \(1630=2\cdot 5\cdot 163\) & b.10 &        \(4001\) & \(111,21,21,111\) &              \(5\) & \(\alpha\) \\
         &         &         &          & \(16300=2^2\cdot 5^2\cdot 163\) & b.10 &        \(0402\) & \(21,111,21,111\) &              \(5\) &            \\
 \(181\) &   \(2\) &   \(5\) & \(1810\) &      \(1810=2\cdot 5\cdot 181\) & b.10 &        \(3010\) & \(111,21,111,21\) &              \(5\) & \(\alpha\) \\
         &         &         &          & \(18100=2^2\cdot 5^2\cdot 181\) & b.10 &        \(0043\) & \(21,21,111,111\) &              \(5\) &            \\
 \(199\) &   \(2\) &   \(5\) & \(1990\) &      \(1990=2\cdot 5\cdot 199\) & b.10 &        \(4001\) & \(111,32,22,111\) &              \(5\) & \(\alpha\) \\
         &         &         &          & \(19900=2^2\cdot 5^2\cdot 199\) & b.10 &        \(4001\) & \(111,32,22,111\) &              \(5\) &            \\
 \(307\) &   \(2\) &   \(5\) & \(3070\) &      \(3070=2\cdot 5\cdot 307\) & b.10 &        \(0043\) & \(21,21,111,111\) &              \(2\) & \(\alpha\) \\
         &         &         &          & \(30700=2^2\cdot 5^2\cdot 307\) & b.10 &        \(4001\) & \(111,21,21,111\) &              \(2\) &            \\
 \(379\) &   \(2\) &   \(5\) & \(3790\) &      \(3790=2\cdot 5\cdot 379\) & b.10 &        \(4001\) & \(111,21,21,111\) &              \(5\) & \(\alpha\) \\
         &         &         &          & \(37900=2^2\cdot 5^2\cdot 379\) & b.10 &        \(4001\) & \(111,21,21,111\) &              \(5\) &            \\
 \(397\) &   \(2\) &   \(5\) & \(3970\) & \(39700=2^2\cdot 5^2\cdot 397\) & b.10 &        \(0043\) & \(21,21,111,111\) &              \(2\) & \(\alpha\) \\
         &         &         &          &      \(3970=2\cdot 5\cdot 397\) & b.10 &        \(0320\) & \(21,111,111,21\) &              \(2\) &            \\
 \(433\) &   \(2\) &   \(5\) & \(4330\) &      \(4330=2\cdot 5\cdot 433\) & b.10 &        \(0320\) & \(21,111,111,21\) &              \(2\) & \(\alpha\) \\
         &         &         &          & \(43300=2^2\cdot 5^2\cdot 433\) & b.10 &        \(0320\) & \(21,111,111,21\) &              \(2\) &            \\
 \(487\) &   \(2\) &   \(5\) & \(4870\) & \(48700=2^2\cdot 5^2\cdot 487\) & b.10 &        \(0043\) & \(22,32,111,111\) &              \(5\) & \(\alpha\) \\
         &         &         &          &      \(4870=2\cdot 5\cdot 487\) & b.10 &        \(3010\) & \(111,32,111,22\) &              \(5\) &            \\
 \(523\) &   \(2\) &   \(5\) & \(5230\) &      \(5230=2\cdot 5\cdot 523\) & b.10 &        \(3010\) & \(111,32,111,21\) &  \(2^2\cdot 5=20\) & \(\alpha\) \\
         &         &         &          & \(52300=2^2\cdot 5^2\cdot 523\) & b.10 &        \(3010\) & \(111,32,111,21\) &  \(2^2\cdot 5=20\) &            \\
 \(541\) &   \(2\) &   \(5\) & \(5410\) &      \(5410=2\cdot 5\cdot 541\) & d.23 &        \(0243\) & \(21,33,111,111\) &  \(2^2\cdot 5=20\) & \(\gamma\) \\
         &         &         &          & \(54100=2^2\cdot 5^2\cdot 541\) & d.23 &        \(3210\) & \(111,33,111,21\) &  \(2^2\cdot 5=20\) &            \\
\hline
\end{tabular}

}

\end{center}
\end{table}


\noindent
The conditions for the \(3\)-class groups of \(k\) and \(\Gamma\)
and the \textit{class number formula}
enforce that all fields \(k\) are of DPF-type \(\beta\)
with \textit{non-radical absolute} principal factors.
The DPF in the tables is always a \textit{cubic residue} modulo the prime \(p\)
(except if it is divisible by \(p\)).


\subsection{Counterexamples of Dedekind species IB}
\label{ss:DedekindIB}

Second, we provide some counterexamples
of Dedekind species IB
with conductors of the shape \(3\Vert f=3pq\),
giving rise to \textit{doublets},
in Table
\ref{tbl:3PQ12} and
\ref{tbl:3PQ15}.
The DPF is always a \textit{cubic residue} modulo the prime \(p\).


\renewcommand{\arraystretch}{1.1}

\begin{table}[ht]
\caption{Doublets, \(m=2\), with \(f=3pq\), \(q\equiv 2\,(\mathrm{mod}\,9)\)}
\label{tbl:3PQ12}
\begin{center}

{\small

\begin{tabular}{|rr|r|ccccc|c|}
\hline
   \(p\) &  \(q\) &    \(f\) &                \(n_i\) &   CT & \(\varkappa_i\) &      \(\alpha_i\) &               DPF &       Type \\
\hline
  \(73\) &  \(2\) &  \(438\) &    \(292=2^2\cdot 73\) & b.10 &        \(4001\) & \(111,21,21,111\) &             \(3\) & \(\alpha\) \\
         &        &          &      \(146=2\cdot 73\) & b.10 &        \(0043\) & \(21,21,111,111\) &             \(3\) &            \\
  \(73\) & \(11\) & \(2409\) &  \(8833=11^2\cdot 73\) & b.10 &        \(0043\) & \(21,21,111,111\) &             \(3\) & \(\alpha\) \\
         &        &          &     \(803=11\cdot 73\) & b.10 &        \(3010\) & \(111,21,111,21\) &             \(3\) &            \\
 \(199\) &  \(2\) & \(1194\) &     \(398=2\cdot 199\) & b.10 &        \(0402\) & \(32,111,21,111\) & \(12=2^2\cdot 3\) & \(\alpha\) \\
         &        &          &   \(796=2^2\cdot 199\) & b.10 &        \(0043\) & \(21,32,111,111\) & \(12=2^2\cdot 3\) &            \\
 \(271\) &  \(2\) & \(1626\) &     \(542=2\cdot 271\) & b.10 &        \(0402\) & \(32,111,22,111\) &             \(3\) & \(\alpha\) \\
         &        &          &  \(1084=2^2\cdot 271\) & b.10 &        \(0402\) & \(32,111,22,111\) &             \(3\) &            \\
 \(487\) &  \(2\) & \(2922\) &  \(1948=2^2\cdot 487\) & b.10 &        \(0043\) & \(21,32,111,111\) & \(12=2^2\cdot 3\) & \(\alpha\) \\
         &        &          &     \(974=2\cdot 487\) & b.10 &        \(0043\) & \(21,32,111,111\) & \(12=2^2\cdot 3\) &            \\
 \(523\) &  \(2\) & \(3138\) &    \(1046=2\cdot 523\) & b.10 &        \(3010\) & \(111,21,111,21\) &             \(3\) & \(\alpha\) \\
         &        &          &  \(2092=2^2\cdot 523\) & b.10 &        \(0043\) & \(21,21,111,111\) &             \(3\) &            \\
 \(541\) &  \(2\) & \(3246\) &  \(2164=2^2\cdot 541\) & d.23 &        \(0243\) & \(21,33,111,111\) &    \(6=3\cdot 2\) & \(\gamma\) \\  
         &        &          &    \(1082=2\cdot 541\) & d.23 &        \(3210\) & \(111,33,111,21\) &    \(6=3\cdot 2\) &            \\
 \(577\) &  \(2\) & \(3462\) &    \(1154=2\cdot 577\) & b.10 &        \(3010\) & \(111,21,111,21\) &             \(3\) & \(\alpha\) \\
         &        &          &  \(2308=2^2\cdot 577\) & b.10 &        \(0043\) & \(21,21,111,111\) &             \(3\) &            \\
 \(613\) &  \(2\) & \(3678\) &  \(2452=2^2\cdot 613\) & b.10 &        \(4001\) & \(111,21,21,111\) &             \(3\) & \(\alpha\) \\
         &        &          &    \(1226=2\cdot 613\) & b.10 &        \(4001\) & \(111,21,21,111\) &             \(3\) &            \\
\hline
\end{tabular}

}

\end{center}
\end{table}


\renewcommand{\arraystretch}{1.1}

\begin{table}[ht]
\caption{Doublets, \(m=2\), with \(f=3pq\), \(q\equiv 5\,(\mathrm{mod}\,9)\)}
\label{tbl:3PQ15}
\begin{center}

{\small

\begin{tabular}{|rr|r|ccccc|c|}
\hline
   \(p\) & \(q\) &     \(f\) &                \(n_i\) &   CT & \(\varkappa_i\) &      \(\alpha_i\) &               DPF &       Type \\
\hline
  \(73\) & \(5\) &  \(1095\) &      \(365=5\cdot 73\) & b.10 &        \(0043\) & \(21,21,111,111\) &             \(3\) & \(\alpha\) \\
         &       &           &   \(1825=5^2\cdot 73\) & b.10 &        \(4001\) & \(111,21,21,111\) &             \(3\) &            \\
 \(541\) & \(5\) &  \(8115\) & \(13525=5^2\cdot 541\) & d.23 &        \(0243\) & \(21,33,111,111\) &   \(15=3\cdot 5\) & \(\gamma\) \\
         &       &           &    \(2705=5\cdot 541\) & d.23 &        \(3210\) & \(111,33,111,21\) &   \(15=3\cdot 5\) &            \\
 \(919\) & \(5\) & \(13785\) &    \(4595=5\cdot 919\) & d.23 &        \(3210\) & \(111,33,111,22\) &             \(3\) & \(\gamma\) \\
         &       &           & \(22975=5^2\cdot 919\) & d.23 &        \(3210\) & \(111,33,111,22\) &             \(3\) &            \\
\hline
\end{tabular}

}

\end{center}
\end{table}


\subsection{Counterexamples of Dedekind species IA}
\label{ss:DedekindIA}

Next, we list lots of counterexamples
of Dedekind species IA
with conductors of the shape \(9\mid f=9pq\),
giving rise to \textit{quartets},
in Table
\ref{tbl:9PQ12} and
\ref{tbl:9PQ15}.
A special remark is due to the conductor \(f=3582\):
Here, two members of the quartet have irregular
\(3\)-class groups of rank \(3\),
namely \((9,9,3)\) and \((9,3,3)\),
the former with radicand \(2388=2^2\cdot 3\cdot 199\) and DPF-type \(\alpha\),
the latter with radicand \(3582=2\cdot 3^2\cdot 199\), DPF-type \(\beta\),
and DPF \(12=2^2\cdot 3\).
The DPF is always \textit{cubic residue} modulo the prime \(p\).


\renewcommand{\arraystretch}{1.1}

\begin{table}[ht]
\caption{Quartets, \(m=4\), with \(f=9pq\), \(q\equiv 2\,(\mathrm{mod}\,9)\)}
\label{tbl:9PQ12}
\begin{center}

{\small

\begin{tabular}{|rr|r|ccccc|c|}
\hline
   \(p\) &  \(q\) &    \(f\) &                         \(n_i\) &   CT & \(\varkappa_i\) &      \(\alpha_i\) &               DPF &       Type \\
\hline
  \(73\) &  \(2\) & \(1314\) &     \(1314=2\cdot 3^2\cdot 73\) & b.10 &        \(0320\) & \(21,111,111,21\) &             \(3\) & \(\alpha\) \\
         &        &          &      \(876=2^2\cdot 3\cdot 73\) & b.10 &        \(3010\) & \(111,21,111,21\) &             \(3\) &            \\
         &        &          &   \(2628=2^2\cdot 3^2\cdot 73\) & b.10 &        \(0402\) & \(21,111,21,111\) &             \(3\) &            \\
         &        &          &        \(438=2\cdot 3\cdot 73\) & b.10 &        \(4001\) & \(111,21,21,111\) &             \(3\) &            \\
 \(109\) &  \(2\) & \(1962\) &    \(1308=2^2\cdot 3\cdot 109\) & b.10 &        \(4001\) & \(111,21,21,111\) &             \(2\) & \(\alpha\) \\
         &        &          &    \(1962=2\cdot 3^2\cdot 109\) & b.10 &        \(0043\) & \(21,21,111,111\) &             \(2\) &            \\
         &        &          &       \(654=2\cdot 3\cdot 109\) & b.10 &        \(0043\) & \(21,21,111,111\) &             \(2\) &            \\
         &        &          &  \(3924=2^2\cdot 3^2\cdot 109\) & b.10 &        \(3010\) & \(111,21,111,21\) &             \(2\) &            \\
 \(127\) &  \(2\) & \(2286\) &    \(1524=2^2\cdot 3\cdot 127\) & b.10 &        \(0320\) & \(21,111,111,21\) &             \(2\) & \(\alpha\) \\
         &        &          &    \(2286=2\cdot 3^2\cdot 127\) & b.10 &        \(0320\) & \(21,111,111,21\) &             \(2\) &            \\
         &        &          &       \(762=2\cdot 3\cdot 127\) & b.10 &        \(0043\) & \(21,21,111,111\) &             \(2\) &            \\
         &        &          &  \(4572=2^2\cdot 3^2\cdot 127\) & b.10 &        \(0043\) & \(21,21,111,111\) &             \(2\) &            \\
 \(199\) &  \(2\) & \(3582\) &      \(1194=2\cdot 3\cdot 199\) & b.10 &        \(3010\) & \(111,32,111,21\) & \(12=2^2\cdot 3\) & \(\alpha\) \\
         &        &          &  \(7164=2^2\cdot 3^2\cdot 199\) & b.10 &        \(0402\) & \(32,111,21,111\) & \(12=2^2\cdot 3\) &            \\
 \(271\) &  \(2\) & \(4878\) &    \(3252=2^2\cdot 3\cdot 271\) & b.10 &        \(3010\) & \(111,32,111,22\) &             \(3\) & \(\alpha\) \\
         &        &          &    \(4878=2\cdot 3^2\cdot 271\) & b.10 &        \(4001\) & \(111,32,22,111\) &             \(3\) &            \\
         &        &          &      \(1662=2\cdot 3\cdot 271\) & b.10 &        \(0043\) & \(22,32,111,111\) &             \(3\) &            \\
         &        &          &  \(9756=2^2\cdot 3^2\cdot 271\) & b.10 &        \(4001\) & \(111,32,22,111\) &             \(3\) &            \\
  \(19\) & \(11\) & \(1881\) &    \(1881=3^2\cdot 11\cdot 19\) & b.10 &        \(4001\) & \(111,21,21,111\) &            \(11\) & \(\alpha\) \\
         &        &          &    \(6897=3\cdot 11^2\cdot 19\) & b.10 &        \(3010\) & \(111,21,111,21\) &            \(11\) &            \\
         &        &          &   \(11913=3\cdot 11\cdot 19^2\) & b.10 &        \(0043\) & \(21,21,111,111\) &            \(11\) &            \\
         &        &          &       \(627=3\cdot 11\cdot 19\) & b.10 &        \(3010\) & \(111,21,111,21\) &            \(11\) &            \\
  \(37\) & \(11\) & \(3663\) &   \(13431=3\cdot 11^2\cdot 37\) & b.10 &        \(4001\) & \(111,21,21,111\) &            \(11\) & \(\alpha\) \\
         &        &          &    \(3663=3^2\cdot 11\cdot 37\) & b.10 &        \(4001\) & \(111,21,21,111\) &            \(11\) &            \\
         &        &          &      \(1221=3\cdot 11\cdot 37\) & b.10 &        \(3010\) & \(111,21,111,21\) &            \(11\) &            \\
         &        &          & \(40293=3^2\cdot 11^2\cdot 37\) & b.10 &        \(3010\) & \(111,21,111,21\) &            \(11\) &            \\
\hline
\end{tabular}

}

\end{center}
\end{table}


\renewcommand{\arraystretch}{1.1}

\begin{table}[ht]
\caption{Quartets, \(m=4\), with \(f=9pq\), \(q\equiv 5\,(\mathrm{mod}\,9)\)}
\label{tbl:9PQ15}
\begin{center}

{\small

\begin{tabular}{|rr|r|ccccc|c|}
\hline
   \(p\) & \(q\) &    \(f\) &                         \(n_i\) &   CT & \(\varkappa_i\) &      \(\alpha_i\) &               DPF &       Type \\
\hline
  \(73\) & \(5\) & \(3285\) &  \(16425=3^2\cdot 5^2\cdot 73\) & b.10 &        \(0043\) & \(21,21,111,111\) &             \(3\) & \(\alpha\) \\
         &       &          &       \(1095=3\cdot 5\cdot 73\) & b.10 &        \(4001\) & \(111,21,21,111\) &             \(3\) &            \\
         &       &          &     \(3285=3^2\cdot 5\cdot 73\) & b.10 &        \(4001\) & \(111,21,21,111\) &             \(3\) &            \\
         &       &          &     \(5475=3\cdot 5^2\cdot 73\) & b.10 &        \(4001\) & \(111,21,21,111\) &             \(3\) &            \\
 \(127\) & \(5\) & \(5715\) &      \(1905=3\cdot 5\cdot 127\) & b.10 &        \(0043\) & \(21,21,111,111\) &             \(5\) & \(\alpha\) \\
         &       &          & \(28575=3^2\cdot 5^2\cdot 127\) & b.10 &        \(0402\) & \(21,111,21,111\) &             \(5\) &            \\
         &       &          &    \(9525=3\cdot 5^2\cdot 127\) & b.10 &        \(0320\) & \(21,111,111,21\) &             \(5\) &            \\
         &       &          &    \(5715=3^2\cdot 5\cdot 127\) & b.10 &        \(0402\) & \(21,111,21,111\) &             \(5\) &            \\
\hline
\end{tabular}

}

\end{center}
\end{table}


\subsection{Prototypes of counterexamples}
\label{ss:Prototypes}

Finally, we present the \textit{prototypes}
of counterexamples with minimal conductors \(f\).
When several non-isomorphic fields share a common conductor
they are collected in a \textit{multiplet} (doublet or quartet), according to
\cite[Thm. 2.1, p. 833, Exm. 1, p. 840]{Ma1992},
\cite[Cor. 3.2, p. 2219]{Ma2014b},
\cite[Thm. 2.2, p. 255]{AMITA2020}.
For the prime divisor \(p\equiv 1\,(\mathrm{mod}\,9)\) of the conductor \(f\),
the \textit{differential principal factorization type},
\cite[Thm. 2.1, p. 254]{AMITA2020},
\cite{AAIMT2022},
of \(\mathbb{Q}(\sqrt[3]{p})\),
necessarily different from type \(\beta\),
is checked.
Since \(p=541\) is the smallest prime with type \(\gamma\),
the type is \(\alpha\), for all primes \(p<541\),
which prohibits a fixed point capitulation at the first component, by
\cite[\S\ 3, Proof of Thm. 2, item (3), p. 55]{IsEM2005}. 

We go further than
\cite{IsEM2005}
by determining the group \(G=\mathrm{Gal}(\mathrm{F}_3^\infty(k)/k)\)
of the maximal unramified pro-\(3\)-extension,
i.e. the Hilbert \(3\)-class field tower, of \(k\),
which is always metabelian \(\mathrm{F}_3^\infty(k)=\mathrm{F}_3^2(k)\)
in our numerical examples.
Except mentioned differently,
the TKT is always b.10, \(\varkappa\sim(0320)\),
with associated ATI \(\alpha\sim\lbrack(21)(1^3)(1^3)(21)\rbrack\),
which uniquely identifies \(G\simeq\langle 243,3\rangle\)
with admissible relation rank \(d_2(G)=4\)
\cite[Thm. 5.1, p. 28]{Ma2015}.
Since \(k\) has
\(3\)-class rank \(\varrho=2\),
signature \((r_1,r_2)=(0,3)\),
torsion free Dirichlet unit rank \(r=r_1+r_2-1=2\),
and contains the primitive cube roots of unity,
indicated by \(\theta=1\),
the bounds for the relation rank of \(G\) are given by
\begin{equation}
\label{eqn:Shafarevich}
2=\varrho\le d_2(G)\le\varrho+r+\theta=2+2+1=5,
\end{equation}
according to Shafarevich
\cite{Sh1964}.
The differential principal factorization type of
\(k=\mathbb{Q}(\zeta,\sqrt[3]{n})\)
is always DPF-type \(\beta\).


\begin{example}
\label{exm:3PQ}
Let \(f=3pq\)
with primes \(p\equiv 1\,(\mathrm{mod}\,9)\)
and \(q\equiv 2,5\,(\mathrm{mod}\,9)\)
(species IB).
\begin{itemize}
\item
For \(q\equiv 2\,(\mathrm{mod}\,9)\), we found four prototypes
(see Table
\ref{tbl:3PQ12}): \\
\(f=3\cdot 73\cdot 2=438\), 
a doublet with radicands \(n\in\lbrace 2^2\cdot 73, 2\cdot 73\rbrace\) (ground state).\\
\(f=3\cdot 199\cdot 2=1194\), 
a doublet with radicands \(n\in\lbrace 2\cdot 199, 2^2\cdot 199\rbrace\) and
ATI \(\alpha\sim\lbrack(32)(1^3)(1^3)(21)\rbrack\) (excited state).\\
\(f=3\cdot 271\cdot 2=1626\), 
a doublet with radicands \(n\in\lbrace 2\cdot 271, 2^2\cdot 271\rbrace\) and
ATI \(\alpha\sim\lbrack(32)(1^3)(1^3)(22)\rbrack\) (excited state in a higher coclass graph
\cite{Ma2013},
see \S\ \ref{s:Coclass}).\\
\(f=3\cdot 541\cdot 2=3246\),
a doublet with radicands \(n\in\lbrace 2^2\cdot 541, 2\cdot 541\rbrace\).
This is the unique instance which admits a fixed point capitulation,
since \(\mathbb{Q}(\sqrt[3]{541})\) is of type \(\gamma\), and indeed
both members of the doublet share the TKT d.23, \(\varkappa\sim (1320)\),
with associated ATI \(\alpha\sim\lbrack(33)(1^3)(1^3)(21)\rbrack\)
and thus \(G\simeq\langle 6561,1990\rangle\) with \(d_2(G)=4\).
\item
For \(q\equiv 5\,(\mathrm{mod}\,9)\), there occurred three prototypes
(see Table
\ref{tbl:3PQ15}): \\
\(f=3\cdot 73\cdot 5=1095\),
a doublet with radicands \(n\in\lbrace 5\cdot 73, 5^2\cdot 73\rbrace\) (ground state).\\
\(f=3\cdot 541\cdot 5=8115\),
a doublet with radicands \(n\in\lbrace 5\cdot 541, 5^2\cdot 541\rbrace\).
Since \(\mathbb{Q}(\sqrt[3]{541})\) is of type \(\gamma\),
a fixed point is admissible in the capitulation type, and
both members of the doublet share the TKT d.23, \(\varkappa\sim (1320)\),
with associated ATI \(\alpha\sim\lbrack(33)(1^3)(1^3)(21)\rbrack\),
again leading to \(G\simeq\langle 6561,1990\rangle\)
\cite{MAGMA6561} (ground state).\\
\(f=3\cdot 919\cdot 5=13785\),
a doublet with radicands \(n\in\lbrace 5\cdot 919, 5^2\cdot 919\rbrace\)
with TKT d.23, \(\varkappa\sim (1320)\), since \(919\in\,\)A363717,
ATI \(\alpha\sim\lbrack(33)(1^3)(1^3)(22)\rbrack\),
(excited state).
\end{itemize}
\end{example} 


\begin{example}
\label{exm:9PQ}
Let \(f=9pq\)
with primes \(p\equiv 1\,(\mathrm{mod}\,9)\)
and \(q\equiv 2,5\,(\mathrm{mod}\,9)\) (species IA).
When \(3\) is cubic residue modulo \(p\), then the principal factors are \(3\) and \(9\),
for each component of the quartets.
Again, the TKT is \(\varkappa\sim (0320)\),
and the (logarithmic) abelian quotient invariants are \(\alpha\sim\lbrack(21)(1^3)(1^3)(21)\rbrack\),
for each component of the quartets.
\begin{itemize}
\item
For \(q\equiv 2\,(\mathrm{mod}\,9)\), we have three prototypes
(see Table
\ref{tbl:9PQ12}): \\
\(f=9\cdot 73\cdot 2=1314\),
quartet with radicands \(n\in\lbrace 2\cdot 3\cdot 73, 2^2\cdot 3\cdot 73, 2^2\cdot 3^2\cdot 73, 2\cdot 3^2\cdot 73\rbrace\) (ground state).\\
\(f=9\cdot 199\cdot 2=3582\),
two fields in a quartet with radicands \(n\in\lbrace 2\cdot 3\cdot 199,2^2\cdot 3^2\cdot 199\rbrace\)
and ATI \(\alpha\sim\lbrack(32)(1^3)(1^3)(21)\rbrack\) (excited state).\\
\(f=9\cdot 271\cdot 2=4878\)
quartet with radicands \(n\in\lbrace 2^2\cdot 3\cdot 271,2\cdot 3^2\cdot 271,2\cdot 3\cdot 271,2^2\cdot 3^2\cdot 271\)
and ATI \(\lbrack(32)(1^3)(1^3)(22)\rbrack\) (excited state in a higher coclass graph
\cite{Ma2013}).
\item
For \(q\equiv 5\,(\mathrm{mod}\,9)\), we got a single prototype in the ground state
(see Table
\ref{tbl:9PQ15}): \\
\(f=9\cdot 73\cdot 5=3285\),
quartet with radicands \(n\in\lbrace 3^2\cdot 5^2\cdot 73, 3\cdot 5\cdot 73, 3^2\cdot 5\cdot 73, 3\cdot 5^2\cdot 73\rbrace\).
\end{itemize}
\end{example} 


\begin{example}
\label{exm:PQ1Q2}
Let \(f=pq_1q_2\)
with primes \(p\equiv 1\,(\mathrm{mod}\,9)\)
and \(q_i\equiv 2,5\,(\mathrm{mod}\,9)\) for \(i=1,2\) (species II).
Again, the capitulation type is \(\varkappa\sim (0320)\),
and the (logarithmic) abelian quotient invariants are \(\alpha\sim\lbrack(21)(1^3)(1^3)(21)\rbrack\),
for each component of the doublets.
\begin{itemize}
\item
For \(q_1\equiv q_2\equiv 2\,(\mathrm{mod}\,9)\), there occurred two prototypes
(see Table
\ref{tbl:PQQ122}): \\
\(f=19\cdot 2\cdot 11=418\), 
a doublet with radicands \(n\in\lbrace 2^2\cdot 11\cdot 19, 2\cdot 11\cdot 19\rbrace\) (ground state).
Since \(11\) is cubic residue modulo \(19\),
the principal factors are \(11\) and \(2^2\cdot 19\), respectively \(11\) and \(2\cdot 19\).
See also \S\ \ref{s:Coclass} with respect to the coclass of the group \(G\).\\
\(f=199\cdot 2\cdot 11=4378\),
a doublet with radicands \(n\in\lbrace 2^2\cdot 11\cdot 199,2\cdot 11^2\cdot 199\rbrace\)
and ATI \(\alpha\sim\lbrack(32)(1^3)(1^3)(22)\rbrack\) (excited state in a higher coclass graph
\cite{Ma2013}).
\item
For the mixed case
\(q_1\equiv 2\,(\mathrm{mod}\,9)\), \(q_2\equiv 5\,(\mathrm{mod}\,9)\), we found four prototypes
(see Table
\ref{tbl:PQQ125}): \\
\(f=19\cdot 5\cdot 11=1045\), 
a doublet with radicands \(n\in\lbrace 5\cdot 11\cdot 19, 5^2\cdot 11\cdot 19\rbrace\) (ground state).
Since \(11\) is cubic residue modulo \(19\),
the principal factors are \(11\) and \(5\cdot 19\), respectively \(11\) and \(11^2\).\\
\(f=199\cdot 2\cdot 5=1990\)
a doublet with radicands \(n\in\lbrace 2\cdot 5\cdot 199,2^2\cdot 5^2\cdot 199\rbrace\)
and ATI \(\alpha\sim\lbrack(32)(1^3)(1^3)(22)\rbrack\) (excited state in a higher coclass graph
\cite{Ma2013}).\\
\(f=523\cdot 2\cdot 5=5230\)
a doublet with radicands \(n\in\lbrace 2\cdot 5\cdot 523,2^2\cdot 5^2\cdot 523\rbrace\)
and ATI \(\alpha\sim\lbrack(32)(1^3)(1^3)(21)\rbrack\) (excited state).\\
\(f=541\cdot 2\cdot 5=5410\)
a doublet with radicands \(n\in\lbrace 2\cdot 5\cdot 541,2^2\cdot 5^2\cdot 541\rbrace\).
Since \(\mathbb{Q}(\sqrt[3]{541})\) is of type \(\gamma\),
a fixed point is admissible in the capitulation type, and
both members of the doublet share the TKT d.23, \(\varkappa\sim (1320)\),
with associated ATI \(\alpha\sim\lbrack(33)(1^3)(1^3)(21)\rbrack\),
leading to \(G\simeq\langle 6561,1990\rangle\)
\cite{MAGMA6561} (ground state).
\item
For \(q_1\equiv q_2\equiv 5\,(\mathrm{mod}\,9)\), there occurred a single prototype
(see Table
\ref{tbl:PQQ155}): \\
\(f=37\cdot 5\cdot 23=4255\), 
a doublet with radicands \(n\in\lbrace 5\cdot 23\cdot 37, 5^2\cdot 23\cdot 37\rbrace\) (ground state).
Since \(23\) is cubic residue modulo \(37\),
the principal factors are \(23\) and \(5\cdot 37\), respectively \(23\) and \(23^2\).
\end{itemize}
\end{example}


\noindent
Although it is not a counter example, we illuminate the impact of
the prime \(p=541\) once more,
with a conductor intentionally skipped from our systematic investigations.
\(541\) is the leading term of the sequence A363717
in the On-line Encyclopedia of Integer Sequences (OEIS)
\cite{OEIS2025}.
The sequence consists of the prime radicands \(p\equiv 1\,(\mathrm{mod}\,9)\)
of pure cubic number fields \(\mathbb{Q}(\sqrt[3]{p})\)
whose normal closure possesses
the differential principal factorization (DPF) type \(\gamma\)
\cite[Thm. 2.1, p. 254]{AMITA2020}.

\begin{example}
\label{exm:9P}
Let \(f=9p\)
with a prime \(p\equiv 1\,(\mathrm{mod}\,9)\).
Except mentioned differently,
the TKT is usually a.1, \(\varkappa\sim(0000)\).
\begin{itemize}
\item
For \(f=9\cdot 19=171\), we have a doublet,
each component with ATI \(\alpha\sim\lbrack(1)(1)(1)(1)\rbrack\),
which uniquely identifies the abelian root \(G\simeq\langle 9,2\rangle\)
with \(d_2(G)=3\)
\cite{BEO2005}.
\item
For \(f=9\cdot 199=1791\), we found a doublet,
each component with ATI \(\alpha\sim\lbrack(21)(11)(11)(11)\rbrack\),
which uniquely identifies mainline \(G\simeq\langle 81,9\rangle\)
with \(d_2(G)=4\).
\item
For \(f=9\cdot 541=4869\), however, there occurred a doublet,
each component with
TKT \(\varkappa\sim (1000)\), having a fixed point,
and ATI \(\alpha\sim\lbrack(33)(11)(11)(11)\rbrack\),
which uniquely identifies \(G\simeq\langle 2187,387\rangle\)
with \(d_2(G)=3\) in depth \(1\)
\cite{BEO2005}.
\end{itemize}
\end{example}


\noindent
However, eventually we point out that
a pure cubic auxiliary field \(\mathbb{Q}(\sqrt[3]{p})\) of DPF-type \(\gamma\)
is \textbf{not a warranty} for \textit{partial} capitulation in the first component
(it is only necessary but not sufficient for the capitulation types
a.2, \(\varkappa\sim (1000)\) and d.23, \(\varkappa\sim (1320)\)).

\begin{example}
\label{exm:PQ}
Let \(f=pq\) with primes
\(p\equiv 1\,(\mathrm{mod}\,9)\) and
\(q\equiv 8\,(\mathrm{mod}\,9)\).
\begin{itemize}
\item
For \(p=1279\) and \(q=17\), the conductor \(f=21743\) gives rise to
a doublet of Dedekind species II,
whose members have \textit{total} capitulation in the first component
with Artin pattern
\(\alpha=\lbrack 32,11,11,11\rbrack\),
\(\varkappa=(0000)\),
although \(1279\) belongs to the sequence A363717.
\end{itemize}
\end{example}


\section{Experimental proof by coclass}
\label{s:Coclass}

\noindent
According to the On-line Encyclopedia of Integer Sequences (OEIS)
\cite{OEIS2025},
the minimal conductors \(f\) of normal closures \(K=\mathbb{Q}(\zeta_3,\sqrt[3]{d})\)
of pure cubic fields \(L=\mathbb{Q}(\sqrt[3]{d})\)
with elementary bicyclic \(3\)-class group \(\mathrm{Cl}_3(K)\simeq (3,3)\)
and assigned coclass \(\mathrm{cc}(M)=0,1,2,3,\ldots\)
of the second \(3\)-class group \(\mathrm{Gal}(\mathrm{F}_3^2(K)/K)\)
are collected in the sequence A380104.

The term for \(\mathrm{cc}(M)=2\) is \(a(3)=418=2\cdot 11\cdot 19\),
which appears in our Table
\ref{tbl:PQQ122}
as a doublet with radicands
\(836=2^2\cdot 11\cdot 19\) and \(4598=2\cdot 11^2\cdot 19\),
both with capitulation type b.10, \(\varkappa=(4001)\),
and logarithmic abelian type invariants \(\alpha=(111,21,\mathbf{21},111)\).
According to
\cite{AoMa2025},
\(\mathrm{cc}(M)=\#(\mathbf{21})-1=3-1=\mathbf{2}\).

The term for \(\mathrm{cc}(M)=3\) is \(a(4)=1626=2\cdot 3\cdot 271\),
which appears in our Table 
\ref{tbl:3PQ12}
as a doublet with radicands
\(542=2\cdot 271\) and  \(1084=2^2\cdot 271\),
both with capitulation type b.10, \(\varkappa=(0402)\),
and logarithmic abelian type invariants \(\alpha=(32,111,\mathbf{22},111)\).
According to
\cite{AoMa2025},
\(\mathrm{cc}(M)=\#(\mathbf{22})-1=4-1=\mathbf{3}\).

This is another disproof of the erroneous results in \cite{IsEM2005},
since the claimed capitulation types 
a.1, \((0000)\), and a.2, \((1000)\),
enforce a second \(3\)-class group \(M=\mathrm{Gal}(\mathrm{F}_3^2(K)/K)\) of maximal nilpotency class,
that is, of coclass \(\mathrm{cc}(M)=1\)
\cite{Ma2013}.
So the erroneous claims implicitly suggest that bigger coclass \(\mathrm{cc}(M)>1\) were impossible,
in contradiction to the terms \(a(3)\) and \(a(4)\) of the sequence A380104.


\tocless\section{Data availability statement}
\label{s:Data}

\noindent
Experimental results communicated in this article and
the source code of Magma program scripts
used for the computations
may be requested from the second author by email. \\


\tocless\section{Acknowledgements}
\label{s:Acknowledgements}

\noindent
The second author acknowledges that his research was supported by
the Austrian Science Fund (FWF): projects J0497-PHY, P26008-N25,
and by the Research Executive Agency of the European Union (EUREA):
project Horizon Europe 2021--2027.\\



\end{document}